\documentclass[12pt]{article} 
\usepackage{amssymb}
\begin{document}\def\ov{\over} \def\t{\tau}
\def\ep{\varepsilon}  \def\bc{\begin{center}} \def\ec{\end{center}}
\def\ee{\end{equation}} \def\a{\alpha}  \def\b{\beta}
\def\be{\begin{equation}} \def\s{\sigma}  \def\n{\emptyset} 
\def\inv{^{-1}} \def\R{\mathbb R} \def\d{\delta} \def\pl{\partial}
\def\R{\mathbb R} \def\S{\mathbb S} \def\ld{\ldots} \def\cd{\cdots}
\def\Z{\mathbb Z} \def\P{\mathbb P} \def\x{\xi} \def\C{\mathcal{C}}
\newcommand{\xs}[1]{\x_{\s(#1)}} \def\B{\mathbb B} \def\noi{\noindent}
\def\si{\s\inv} \newcommand{\abxs}[1]{\x_{|\s|(#1)}} 
\def\BSN{\B_N\backslash\S_N} \def\ph{\varphi} \def\l{\ell}
\def\m{\mu} \def\n{\nu} \def\lb{\Big[} \def\rb{\Big]}
\def\A{\mathbb A} \def\sp{\vspace{1ex}} 
\newcommand{\SX}[2]{S(\x_{#1},\,\x_{#2})}

\hfill May 17, 2012 

\bc{\bf \large The Bose Gas and Asymmetric Simple Exclusion Process\\ \vspace{1ex}
on the Half-Line }\ec

\begin{center}{\bf \large Craig A.~Tracy}\\
{\it Department of Mathematics \\
University of California\\
Davis, CA 95616, USA\\
email:\texttt{tracy@math.ucdavis.edu}}\end{center}

\begin{center}{\bf \large Harold Widom}\\
{\it Department of Mathematics\\
University of California\\
Santa Cruz, CA 95064, USA\\
email: \texttt{widom@ucsc.edu}}\end{center}

\bc{\bf I. Introduction}\ec
\sp

In this paper we solve two closely related equations:  (1) the time-dependent Schr\"odinger equation for a system of one-dimensional bosons interacting
via a delta-function potential with particles confined to the half-line $\R^+$; and (2) the Kolmogorov forward equation (master equation) for the one-dimensional asymmetric simple exclusion process (ASEP) with particles confined
to the nonnegative integers $\Z^+$.  These give explicit formulas for (1) Green's function for the $\d$-function gas problem on $\R^+$; and (2) the transition probability for the half-line ASEP. These are both for systems with a finite number of particles.

We use coordinate Bethe Ansatz appropriately modified to account for confinement of the particles to the half-line.  Formulas derived previously for the line \cite{tw1,tw2} are sums over the permutation group $\S_N$. We shall find that in our formulas for the half-line the group $\S_N$ (Weyl group $\A_{N-1}$) gets replaced by the Weyl group $\B_N$. The fact that in the Bethe Ansatz for particles interacting via a $\d$-function potential on the half-line the group  $\A_{N-1}$ is replaced by the group $\B_N$  goes back to Gaudin \cite{Ga}. See also \cite{GS, HO} for further developments.

The Lieb-Liniger \cite{LL} $\d$-function gas model and the ASEP have
recently attracted  much  attention due to the close relationship of these models to exact solutions of the Kardar-Parisi-Zhang (KPZ) equation.  Regarding these connections with KPZ, we refer the reader to
two recent reviews \cite{Co, Do} and references therein.

Here is a description of the results. We first state the earlier results for $\R$ and $\Z$ and then indicate the analogues for $\R^+$ and $\Z^+$. The detailed statements and proofs are in Sections II and III, respectively.

For the $\d$-function Bose gas on $\R$, the Lieb-Liniger Hamiltonian for an $N$-particle system is
\[ H=-\sum_{j=1}^N{\pl^2\ov\pl x_j^2} + 2c\sum_{j<k} \d(x_j-x_k). \]
One seeks solutions to the time-dependent Schr\"odinger equation
\be H\Psi = i\,{\pl\Psi\ov\pl t} \label{eq1}\ee
that satisfies the initial condition
\be \Psi(x_1,\ldots,x_N;0)=\prod_j\d(x_j-y_j)\label{ic1}\ee
where $y_1<\cdots<y_N$. (From these we can produce the solution with arbitrary initial condition.) We consider the repulsive case $c>0$. 

The solutions are to be symmetric in the coordinates $x_j$ (the Bose condition). From this one sees that it is enough to solve (\ref{eq1}) in the region $x_1<\cdots<x_N$,
subject to the boundary conditions
\be \left({\pl\ov\pl x_{j+1}}-{\pl\ov\pl x_j}\right)\Psi\left\vert_{x_{j+1}=x_j}\right. = c\Psi\left\vert_{x_{j+1}=x_j}\right..\label{bc1} \ee

The solution was found in \cite{tw2}, and was given as the sum over the permutation group $\S_N$ of multiple integrals. Specifically, define
\be S(k)=-{c-ik\ov c+ik},\ \ \  \ep(k)=k^2.\label{S1}\ee
For $\s\in\S_N$ an {\it inversion} in $\s$ is an ordered pair
$(\s(i),\,\s(j))$ in which $i<j$ and $\s(i)>\s(j)$.
We set
\be A_\s=\prod\left\{ S(k_a-k_b): (a,\,b)\>\> \textrm{is an inversion in}\>\> \s\right\}.\label{A1}\ee
It was shown that $\Psi_y(x;t)$, the solution of (\ref{eq1}) that satisfies the initial condition
(\ref{ic1}) and boundary conditions (\ref{bc1}), is 
\be \Psi_y(x;t)=\sum_{\s\in\S_N}\int_{\R}\cdots\int_{\R} A_\s\, \prod_{j=1}^N e^{i k_{\s(j)} x_j} \, \prod_{j=1}^Ne^{-i k_j y_j-it\sum_j \ep(k_j)}\,
dk_1\cdots dk_N.\footnote{All integrals over $\R$ are given the factor $(2\pi)\inv$.} \label{Psi}\ee

In the asymmetric simple exclusion process, particles are at integer sites on the line. Each particle waits exponential time, and then with probability $p$ it moves one step to the right if the site is unoccupied, otherwise it does not move; and with probability $q=1-p$ it moves one step to the left if the site is unoccupied, otherwise it does not move. For $N$-particle ASEP a possible configuraion is given by 
\[X=\{x_1,\ld,x_N\},\quad x_1<\cd<x_N,\quad (x_i\in\Z).\]
The $x_i$ are the occupied sites. Denote by $\P_Y(X;t)$ the probability that at time~$t$ the system is in configuration~$X$ given that initially it was in configuration 
\[Y=\{y_1,\ld,y_N\}.\]
The probability $\P_Y(X;\,t)$ is the solution of the differential equation 
\[{d\ov dt}\,u(X;t)=\sum_{i=1}^N
\lb p\,u(x_i-1)\,(1-\d(x_i-x_{i-1}-1))
+q\,u(x_i+1)\,(1-\d(x_{i+1}-x_i-1))\]
\be -p\,u(x_i)\,(1-\d(x_{i+1}-x_i-1))-q\,u(x_i)\,(1-\d(x_i-x_{i-1}-1))\rb\label{eq2}\ee
that satisfies the initial condition
\be u(X;\,0)=\d_Y(X).\label{ic2}\ee
(In the $i$th summand in (\ref{eq2}) entry $i$ is displayed and entry $j$ is $x_j$ when $j\ne i$. Any $\d$-term involving $x_0$ or $x_{N+1}$ is replaced by zero.)

Equation (\ref{eq2}) holds if $u$ satisfies 
\be{d\ov dt}\,u(X;t)
=\sum_{i=1}^N\,[p\,u(x_i-1)+q\,u(x_i+1)-u(x_i)]\label{eq3}\ee
for all $x_1,\ld,x_N$, and the boundary conditions 
\be p\,u(x_i,x_i)+q\,u(x_i+1,x_i+1)-u(x_i,x_i+1)=0 \label{bc2}\ee
for $i=1,\ld,N-1$. (Here entries $i$ and $i+1$ are displayed.) One sees this by subtracting the right sides of (\ref{eq2}) and (\ref{eq3}).

The solution was found in \cite{tw1,tw3}, and was given as the sum over $\S_N$ of multiple integrals. We now define
\be S(\x,\,\x')=-{p+q\,\x\,\x'-\x\ov p+q\,\x\,\x'-\x'},\ \ \  \ep(\x)=p\,\x\inv+q\,\x-1,\label{S2}\ee
\be A_\s=\prod\left\{ S(\x_a,\,\x_b)): (a,\,b)\>\> \textrm{is an inversion in}\>\> \s\right\}.\label{A2}\ee
The result \cite[Theorem 2.1]{tw1} was that if $q\ne0$ then
\be \P_Y(X;t)=\sum_{\s\in\S_N} \int_{\C{}^N} A_\s(\x)\,\prod_i\xs{i}^{x_i}\;\prod_i\Big(\x_i^{y_i-1}e^{\ep(\x_i)\,t}\Big)\,d\x_1\cd d\x_N.\footnote{All integrals over $\C$ are given the factor $(2\pi i)\inv$.} \label{PXY}\ee
Here $\C$ is a circle about zero which is so large that the $S$-factors are analytic for the $\x_i$ inside and on $\C$.\footnote{The assumptions were actually that $p\ne0$ and that $\C$ was so small that the $S$-factors were analytic for the $\x_i$ inside and on $\C$. That it also holds for $q\ne0$ and $\C$ large was explained in the remark after \cite[Lemma 2.4]{tw1}.}  

We shall show that for the $\d$-function Bose gas and ASEP on the half-line, both suitably interpreted, there are fomulas analogous to (\ref{Psi}) and (\ref{PXY}). Instead of the sums being taken over $\S_N$ they are taken over the Weyl group $\B_N$, which is most conveniently interpreted here as the group of {\it signed permutations}. These are functions $\s:[1,\,N]\to[-N,\,-1]\cup[1,\,N]$ such that $|\s|$ is a permutation in the usual sense.\footnote{These can be identified in an obvious way with bijections $\s:[-N,\,N]\to[-N,\,N]$ satisfying $\s(-i)=-\s(i)$. From this the group structure of $\B_N$ is clear.} 

An {\it inversion} in $\B_N$ is defined to be a pair $(\pm\s(i),\s(j))$ 
with $i<j$ such that $\pm\s(i)>\s(j)$. For example, if $\s=(-3,\, 1,\, -2)$ then the inversions are $(3,\,1),\,(3,\,-2),\linebreak
(-1,\,-2)$, and $(1,\,-2)$. 

In our final formulas, for the analogues of (\ref{A1}) and (\ref{A2}) we will define $k_a=-k_{-a}$ and $\x_a=\t/\x_{-a}$ when $a<0$. (Here $\t=p/q$). There will be also be extra factors required, to take account of the behavior at zero. They are given in formulas (\ref{AR+}) and (\ref{AZ+}) below. Just as the proof for the Bose gas on $\R$ was more straightforward than for ASEP on $\Z$, so will the proof for the Bose gas on $\R^+$ be more straightforward than for ASEP on $\Z^+$.

\bc{\bf II. Bose Gas on the Half-Line}\ec

Here in the equation (\ref{eq1}), the initial condition (\ref{ic1}), and the boundary conditions (\ref{bc1}) we assume that $0<x_1<\cd<x_N$, and in the initial condition we assume that $0<y_1<\cd<y_N$.
The half-line restriction includes the  {\it hard-wall boundary condition} 
\be \Psi(0+,\,x_2,\ldots,x_N;t)=0 \label{bc3}.\ee

For this half-line case we define
\be A_\s:=(-1)^{\#\{i:\,\s(i)<0\}}\times\prod_{{\rm inversions}\ (a,b)}S(k_{a}-k_{b}).\label{AR+}\ee
Recall that the inversions are those in $\B_N$ and we define $k_a=-k_{-a}$ when $a<0$.

\noi{\bf Theorem}. The solution of equation (\ref{eq1}), with initial condition (\ref{ic1}) and boundary conditions (\ref{bc1}) and (\ref{bc3}), is
\be\Psi_y(x;t)=\sum_{\s\in\B_N}\int_{\R}\cdots\int_{\R} A_\s\, \prod_{j=1}^N e^{i k_{\s(j)} x_j} \,\prod_{j=1}^N e^{-ik_j y_j-i t \ep(k_j)}\,
dk_1\cdots dk_N.\label{Psi1}\ee

\noi{\bf Proof}. It is easy to see that (\ref{Psi1}) satisfies the equation (\ref{eq1}). This is true no matter how the $A_\s$ are defined. Then, as in \cite{tw2}, the boundary conditions (\ref{bc1}) will be satisfied if for all $\s\in\B_N$,
\be{A_{T_i\s}\ov A_\s}=S(k_{\s(i+1)}-k_{\s(i)}),\label{As1}\ee
where $T_i$ interchanges the values of $\s(i)$ and $\s(i+1)$.

We first show that these are satisfied by 
\[A_\s^0:=\prod_{{\rm inversions}\ (a,b)}S(k_{a}-k_{b}).\]
The only things that can change when we apply $T_i$ to $\s$ are the inversions due to a pair $(\s(i),\,\s(i+1))$. In the following table we list all possibilites for such pairs, then the ordering of the absolute values of the constituents, then the inversions they give rise to, then the inversions for the pair we get after switching the entries (i.e., applying~$T_i$), then the product of the $S$-factors coming from the second set of inversions divided by the product of the $S$-factors coming from the first set of inversions. 
table, $a,\,b>0$.
\[\begin{array}{cclll}{\rm pair}&{\rm ordering}&{\rm inversions}&\textrm{inversions after switch}\ \ \ &\textrm{ratio of $S$-factors}\\
&&&&\\
a,\,b&a<b&&(b,a)&S(k_b-k_a)\\
a,\,b&a>b&(a,b)& &S(k_a-k_b)\inv\\
-a,\,b&a<b& &(b,-a)&S(k_b+k_a)\\
-a,\,b&a>b&(a,b)&(b,-a),\,(-b,-a)&S(k_b+k_a)\\
a,\,-b&a<b&(a,-b),\,(-a,-b)\ \ \ &(b,a)&S(k_a+k_b)\inv\\
a,\,-b&a>b&(a,-b)& &S(k_a+k_b)\inv\\
-a,\,-b&a<b&(-a,-b),\,(a,-b)&(b,-a)&S(k_b-k_a)\inv\\
-a,\,-b&a>b&(a,-b)&(-b,-a),\,(b,-a)&S(k_a-k_b).\\
\end{array}\]
If we use $S(k)\inv=S(-k)$ we deduce that (\ref{As1}) is satisfied by the $A_\s^0$ in all cases.

Relations (\ref{As1}) also hold for the $A_\s$ since they hold for the $A_\s^0$ and the $T_i$ have no effect on the first factor in (\ref{AR+}). Thus we have (\ref{bc1}).

Now we show that (\ref{bc3}) also holds. Pair $\s$ and $\s'$ if $\s'(1)=-\s(1)$ and $\s'(i)=\s(i)$ for $i>1$. The inversions for $\s$ and $\s'$ are the same, so $A_\s^0=
A_{\s'}^0$. Since the number of negative numbers in the ranges of $\s$ and $\s'$ differ by one, we have $A_{\s}+A_{\s'}=0$ for each pair $(\s,\,\s')$. 
When $x_1=0$ the remaining factors in the integrands in (\ref{Psi}) are the same for $\s$ and $\s'$, so the sum of the two integrands equals zero. Thus (\ref{bc3}) holds.

It remains to verify the initial condition (\ref{ic1}). It is enough to show that 
\[\int_{\R}\cdots\int_{\R} A_\s\, \prod_{j=1}^N e^{i k_{\s(j)} x_j}\ \prod_{a=1}^N e^{-i k_a y_a} \,dk_1\cdots dk_N =0\]
when $\s$ is not the identity permutation. We know from \cite[Sec. 2.2]{tw2} that this is true when all $\s(i)>0$, so we may assume this is not the case. Let $-b=\s(j)$ be the most negative value of $\s$. Then the only inversions involving $\pm b$ are of the form $(a,\,-b)$ or $(b,\,a)$. Therefore all $S$-factors involving $k_b$ are of the form $S(k_b\pm k_a)$. The other factors involving $k_b$ combine as $e^{-ik_b\,(x_j+y_b)}$. Since $S$ is analytic in the lower half-plane and $x_j,\,y_b>0$, integration with respect to $k_b$ gives zero. Thus, (\ref{ic1}) is satisfied.

\bc{\bf III. ASEP on the Half-Line}\ec

Here particles are restricted to the nonnegative integers $\overline{\Z^+}$, where probabilities for the sites $x\in\Z^+$ are as before while a particle at 0 hops to the right with probability $p$ (if site 1 is unoccupied) and stays at 0 with probability $q$. We assume $q\ne0$.

Now in the $i=0$ term on the right side of (\ref{eq2}) the first and last summands have the factor $1-\d(x_1)$, because a particle cannot move left from 0. This equation would hold if $u$ satisfied (\ref{eq3}) and (\ref{bc2}) as before, and in addition the boundary condition
\be u(0,x_2,\ld,x_N)-\t\,u(-1,x_2,\ld,x_N)=0.\label{bc4}\ee

With $S(\x,\,\x')$ defined by (\ref{S2}), we expect an analogue of (\ref{PXY}) where the sum is over $\B_N$ and the analogue of the $A_\s$ would contain factors $S(\x_a,\,\x_b)$ with $(a,\,b)$ an inversion in~$\s$. Recall that we use
\be\x_{-a}=\t/\x_{a},\label{ximinus}\ee
where $\t=p/q$. But there is a difficulty: since
\[S(\x,\,\t/\x')={\x+\x'-p\inv\,\x\x'\ov\x+\x'-q\inv},\]
a factor $S(\x_a,\,\x_b)$ has singularities on the contours when $a>0$ and $b<0$ since $\C\cap(q\inv-\C)\ne\emptyset$.\footnote{Inside small contours, which we may take in (\ref{PXY}), there will be no singularities. But if we take small contours in (\ref{PXY1}) below the initial condition (\ref{ic2}) will no longer be satisfied, even when $N=1$, as is easily seen from (\ref{N=1}). That explains the large contours.} 
This problem is avoided if the $\x_a$ run over circles with center $1/2q$ and different radii. However, the argument that follows requires that the domain of integration be symmetric in the $\x_a$. Therefore we average over all choices of radii for the~$\x_a$. 
To be precise, fix $R_1<\cd<R_N$ with the $R_a$ large, and denote by $\C_a$ the circle with center $1/2q$ and radius $R_a$. We take as our domain of integration 
\be\bigcup_{\m\in\S_N}\,\C_{\m(1)}\times\cd\times\C_{\m(N)}.\label{domain}\ee

We now define 
\be A_\s=\prod_{\s(i)<0}r(\xs{i})\,\times\,\prod_{{\rm inversions}\ (a,b)}S(\x_{a},\,\x_{b}),\label{AZ+}\ee
where
\[r(\x):=-{1-\x\ov1-\t\,\x\inv}.\]

\noi{\bf Theorem}. For ASEP on the half-line we have
\be \P_Y(X;t)={1\ov N!}\sum_{\s\in\B_N} \int\cd\int A_\s(\x)\,\prod_i\xs{i}^{x_i}\;\prod_i\Big(\x_i^{y_i-1}e^{\ep(\x_i)\,t}\Big)\,d\x_1\cd d\x_N,\label{PXY1}\ee
where $A_\s$ is given by (\ref{AZ+}) and the domain of integration by 
(\ref{domain}).
\sp

\noi {\bf Remark}. The ASEP on the half-line for $N=1$ is one of the simplest examples of a birth-and-death process.\footnote{See, e.g., \cite[Chapter 17]{Fe}.} Formula (\ref{PXY1}) for the transition probability then,
\be\P_y(x;t)=\int_\C\left[ \x^{x-y-1} -\left({1-\t/\x\ov 1-\x}\right)\,\t^x\,\x^{-x-y-1}\right]\,e^{\ep(\x)t}\, d\x, \label{N=1}\ee
is a known one \cite[Chapter 4]{LR} for this special case of the birth-and-death process where the transition rates depend upon the states.
\sp

\noi{\bf Proof}. We must verify equation (\ref{eq3}), boundary conditions (\ref{bc2}) and (\ref{bc4}), and initial condition (\ref{ic2}). 

Each term on the right side of (\ref{PXY1}) satisfies (\ref{eq3}) no matter what the $A_\s$ are,\footnote{This uses the fact that $\ep(\x)$ is invariant under the mapping $\x\to\t/\x$. In the preceding section we used the fact that $\ep(k)$ was invariant under $k\to-k$.} so (\ref{eq3}) holds for the sum.

As in \cite{tw1}, boundary conditions (\ref{bc2}) will be satisfied if the $A_\s$ satisfy
\be{A_{T_i\s}\ov A_\s}=S(\xs{i+1},\,\xs{i}).\label{As2}\ee

We first show these are satisfied by
\[A_\s^0:=\prod_{{\rm inversions}\ (a,b)}S(\x_{a},\,\x_{b}).\]
Using the easily verified identity 
\be S(\x_{-a},\,\x_{-b})=S(\x_b,\,\x_a),\label{Sab}\ee
we find that the analogue of the table in the last section is here
 
\[\begin{array}{cclll}{\rm pair}&{\rm ordering}&{\rm inversions}&\textrm{inversions after switch}\ \ \ &\textrm{ratio of $S$-factors}\\
&&&&\\
a,\,b&a<b&&(b,a)&S(\x_b,\,\x_a)\\
a,\,b&a>b&(a,b)& &S(\x_a,\,\x_b)\inv\\
-a,\,b&a<b& &(b,-a)&S(\x_b,\,\x_{-a})\\
-a,\,b&a>b&(a,b)&(b,-a),\,(-b,-a)&S(\x_b,\,\x_{-a})\\
a,\,-b&a<b&(a,-b),\,(-a,-b)\ \ \ &(b,a)&S(\x_a,\,\x_{-b})\inv\\
a,\,-b&a>b&(a,-b)& &S(\x_a,\,\x_{-b})\inv\\
-a,\,-b&a<b&(-a,-b),\,(a,-b)&(b,-a)&S(\x_{-a},\,\x_{-b})\inv\\
-a,\,-b&a>b&(a,-b)&(-b,-a),\,(b,-a)&S(\x_{-b},\,\x_{-a}).\\
\end{array}\]
If we use 
\be S(\x,\,\x')\,S(\x',\,\x)=1\label{Sprod}\ee
we deduce that (\ref{As2}) is satisfied by the $A_\s^0$ in all cases.

The other boundary condition, (\ref{bc4}), will be satisfied if
\[\sum_{\s\in\B_N}A_\s \,(1-\t\,\xs{1}\inv)=0.\]
By (\ref{ximinus}) this may be written
\[\sum_{\s\in\B_N}A_\s \,(1-\x_{-\s(1)})=0.\] 
This will hold if, with $\s$ and $\s'$ paired as in the last section so that $\s'(1)=-\s(1)$, we have
\[{A_{\s'}\ov1-\x_{\s'(1)}}=-{A_\s\ov 1-\xs{1}}.\]
We find that if we define $A_\s$ as we did in (\ref{AZ+})
then this relation holds. (Observe that the product of $S$-factors for $\s$ and $\s'$ are the same since they have the same inversions.) And (\ref{As2}) continues to hold since the first factor in $A_\s$ is not affected by the $T_i$.

We now know that the equation and boundary conditions are satisfied, and it remains to verify the initial condition (\ref{ic2}). 

Denote by $I(\s)$ the $\s$-summand of the right side of (\ref{PXY1}) with $t=0$. What we have to show is that
\[\sum_{\s\in\B_N}I(\s)=\d_Y(X).\]
We know from \cite{tw3} that the sum over $\s\in\S_N$ equals $\d_Y(X)$. Therefore it remains to show that
\be\sum_{\s\in\BSN}I(\s)=0.\label{sum}\ee
Henceforth we consider only $\s\in\B_N\backslash\S_N$, those for which some $\s(j)<0$.

In the proof of this we use $i$ when $\s(i)>0$ and $j$ when $\s(j)<0$. It is convenient to make the substitutions $\x_a\to\x_a+1/2q$, so that the $\C_a$ become circles with center zero and the interesting part of the integrand becomes
\be\prod_{i}\xs{i}^{x_i-y_{\s(i)}-1}\,
\prod_{j}\abxs{j}^{-x_j-y_{|\s|(j)}-2}\,
\prod_{i<j}{\xs{i}\,\abxs{j}\ov\xs{i}+\abxs{j}}\,
\prod_{j<j'}{\abxs{j}\,\abxs{j'}\ov\abxs{j}+\abxs{j'}}.\label{prod}\ee
By ``the interesting part'' we mean enough to show the location of the poles when we expand contours, and the orders of magnitude of the factors at infinity. (Observe that each $r(\xs{j})$ is $O(\abxs{j}\inv)$ at infinity, which accounts for the exponents $-2$ instead of $-1$ in the second product. It follows from identity (\ref{Sab}) that an $S$-factor $S(\x_a,\,\x_b)$ coming from an inversion $(a,\,b)$ has poles only when $a>0,\,b<0$, which accounts for the last two products.)

We shall integrate with respect to some of the variables by expanding their contours, leaving us with lower-order integrals, the residues at the poles that are passed. After two steps we will be left with subintegrals in which two of the variables are equal. 

Here is how we do it. Keep in mind our use of the indices $i$ and $j$. If $j_0$ is the largest $j$ we integrate with respect to $\abxs{j_0}$ by expanding its contour. We pass simple poles at the $-\xs{i_1}$ with $i_1<j_0$ and at the $-\abxs{j_1}$ with 
$j_1<j_0$.\footnote{Not all of these cases need actually arise. For example if $\s$ takes only one negative value then the only poles passed in this integration are at the $-\xs{i_1}$ with $i_1<j_0$. Alao, for each $i_1$ or $j_1$ only some of the contours in (\ref{domain}) contribute. See the discussion after the end of the proof. } After the integration the exponent of $\xs{i_1}$ is reduced by $x_{j_0}+y_{|\s|(j_0)}$ and so the resulting exponent is at most $-2$, since $i_1<j_0$. The resulting exponent of $\abxs{j_1}$ is even more negative. 

The residue at $\abxs{j_0}=-\xs{i_1}$ we integrate with respect to $\xs{i_1}$ by expanding its contour. We pass simple poles at some $-\abxs{j_1}$ with $j_1\ne j_0$, at some $\xs{i_2}$ with $i_2\ne i_1$, and at some $\abxs{j_1}$ with $j_1\ne j_0$. Since after the first integration we have $\abxs{j_0}=-\xs{i_1}$, after this second integration we have $\abxs{j_1}=\abxs{j_0}$ in the first case, 
$\xs{i_1}=\xs{i_2}$ in the second case, and $\abxs{j_1}=\xs{i_1}$ in the third case.

The residue at $\abxs{j_0}=-\xs{j_1}$ we integrate with respect to $\abxs{j_1}$. We pass poles at some $-\xs{i_1}$, at some $\xs{i_1}$, at some $-\abxs{j_2}$ with $j_2\ne j_0,\,j_1$, and at some $\abxs{j_2}$ with $j_2\ne j_0,\,j_1$. Since after the first integration we have $\abxs{j_0}=-\xs{j_1}$, after the second integration we have $\abxs{j_0}=\xs{i_1}$ in the first case, $\abxs{j_1}=\xs{i_1}$ in the second case, $\abxs{j_0}=\abxs{j_2}$ in the third case, and $\abxs{j_1}=\abxs{j_2}$ in the last case.

Thus after two integrations any nonzero $I(\s)$ is represented as the sum of subintegrals in each of which two of the variables are equal. If we take two different~$\s$ we get different subintegrals. 

Consider a subintegral where, say, $\abxs{j_1}=\xs{i_1}$, and correspondingly the subintegral where $\x_{|\s'|(j_1)}=\x_{\s'(i_1)}$. Both the domains of integration and the integrands are different for the two. But consider the permutations
\be\s=(1,\ -2,\ 3,\ 5,\ -4),\ \ \ \s'=(1,\ -5,\ 3,\ 2,\ -4).\label{example}\ee
All entries of $\s$ and $\s'$ are the same except for entries 2 and 4. With $j_1=2$ and $i_1=4$ we have $\s(j_1)=-2,\ \s(i_1)=5$ and $\s'(j_1)=-5,\ \s(i_1)=2$. In both subintegrals $\x_2=\x_5$. If we interchange the variables $\x_2$ and $\x_5$ in the $\s'$-integral then the domains of integration become the same for the two, by the symmetry of the original domain of integration,\footnote{We explain this after the end of the proof.} and the integrands themselves become almost the same; what is different are only the $S$-factors arising from the two permutations.

This is quite general. For given $a,\,b>0$, say that $\s$ and $\s'$ are $(a,\,b)$-paired if they agree except for the positions of $\pm a$ and $\pm b$, and the positive numbers $a$ and~$b$ are interchanged. The permuatations in (\ref{example}) are $(2,\,5)$-paired. For subintegrals in which the variables $\x_a$ and $\x_b$ are equal, if permutations $\s$ and $\s'$ are $(a,\,b)$-paired then what is different in the integrands after interchanging the variables $\x_a$ and $\x_b$ in the $\s'$-integral are only the $S$-factors arising from the two permutations. We shall show that when $\x_a=\x_b$ the product of $S$-factors for $\s$ and $\s'$ are negatives of each other. (This whether or not we interchange the variables.) Thus, the sum of the two subintegrals equals zero. 

In what follows we always assume that $a<b$ and that $\pm a$ appears before $\pm b$ in~$\s$. (Otherwise we reverse the roles of $\s$ and $\s'$.) There are four situations, which depend on which pair of signs occurs. We consider them separately and, with obvious notation, we denote the four cases by $(+,\,+),\ (-,\,-),\ (+,\,-),\ (-,\,+)$. In example (\ref{example}) the sign pair is $(-,\,+)$.

In each case we show that when $\x_a=\x_b$ the products of $S$-factors involving only $\pm a$ and $\pm b$ are negatives of each other, and that for any $c\ne \pm a,\,\pm b$ the products of $S$-factors involving $\x_{\pm c}$ and $\x_{\pm a}$ or $\x_{\pm b}$ are equal. This will give the desired result. If $\si(c)$ is outside the interval $(\si(\pm a),\,\si(\pm b))$ the $S$-factors in question are the same for $\s$ and $\s'$, so we will always assume $\si(c)$ is inside this interval. For each of the four cases there will be five subcases, depending on the position of $c$ relative to $\pm a$ and $\pm b$, with the results displayed in tables. The first column will tell where $c$ is relative to $\pm a$ and $\pm b$, the second column will give the product of $S$-factors involving $\x_{\pm c}$ and either $\x_{\pm a}$ or $\x_{\pm b}$ for $\s$, and the fourth column will give the corresponding product for $\s'$. 
\sp

\noi{\bf The case \boldmath$(+,\,+)$}: The only $S$-factor involving only $\pm a$ and $\pm b$ is $S(\x_b,\,\x_a)$ for $\s'$. This equals $-1$ when $\x_a=\x_b$. For $c\ne \pm a,\,\pm b$ the table described above is 

\[\begin{array}{lllll}
c<-b && \SX{a}{c}\,\SX{-a}{c}\,\SX{-c}{b} &&  \SX{b}{c}\,\SX{-b}{c}\,\SX{-c}{a}\\
-b<c<-a && \SX{a}{c}\,\SX{-a}{c} &&  \SX{b}{c}\,\SX{-c}{a}\\
-a<c<a && \SX{a}{c} &&  \SX{b}{c}\\
a<c<b && 1 &&  \SX{b}{c}\,\SX{c}{a}\\
c>b && \SX{c}{b} &&  \SX{c}{a}\\
\end{array}\]
If we use identity (\ref{Sab}) in the second row and identity (\ref{Sprod}) in the fourth, we see that the two columns are the same when $\x_a=\x_b$. 
\sp

\noi{\bf The case \boldmath$(-,\,-)$}: Now the $S$-factors involving only $\pm a$ and $\pm b$ are  $\SX{a}{-b}\,\SX{-a}{-b}$ for $\s$ and 
$\SX{b}{-a}$ for $\s'$. These are clearly negatives when $\x_a=\x_b$. For 
$c\ne\pm a,\,\pm b$ we have the table
\[\begin{array}{lllll}
c<-b && \SX{a}{c}\,\SX{-a}{c}\,\SX{-c}{b} && \SX{b}{c}\,\SX{-b}{c}\,\SX{-c}{a} \\
-b<c<-a &&  \SX{a}{c}\,\SX{-c}{-b}\,\SX{-a}{c}\,\SX{c}{-b}&& \SX{b}{c}\,\SX{-c}{-a} \\
-a<c<a && \SX{a}{c}\,\SX{c}{-b}\,\SX{-c}{-b}&& \SX{b}{c}\,\SX{c}{-a}\,\SX{-c}{-a}\\
a<c<b && \SX{-c}{-b}\, \SX{c}{-b}&& \SX{b}{c}\,\SX{c}{-a} \\
c>b  &&  \SX{c}{-b}&& \SX{c}{-a}\\
\end{array}\]
Using (\ref{Sab}) and (\ref{Sprod}) we see that the two columns are the same when $\x_a=\x_b$.
\sp

\noi{\bf The case \boldmath$(+,\,-)$}: The $S$-factors involving only $\pm a$ and $\pm b$ are $\SX{a}{-b}\,\SX{-a}{-b}$ for $\s$ and 
$\SX{b}{-a}$ for $\s'$, which are again clearly negatives of each other when $\x_a=\x_b$. For $c\ne\pm a,\,\pm b$ we have the table
\[\begin{array}{lllll}
c<-b && \SX{a}{c}\,\SX{-a}{c}\,\SX{-c}{-b} && \SX{b}{c}\,\SX{-b}{c}\,\SX{-c}{-a} \\
-b<c<-a && \SX{a}{c}\,\SX{-a}{c}\,\SX{c}{-b}&& \SX{-c}{-a} \\
-a<c<a && \SX{a}{c}\,\SX{c}{-b}\,\SX{-c}{-b}&& \SX{b}{c}\,\SX{c}{-a}\,\SX{-c}{-a}\\
a<c<b && \SX{-c}{-b}\, \SX{c}{-b}&& \SX{b}{c}\,\SX{c}{-a} \\
c>b  &&  \SX{c}{-b}&& \SX{c}{-a}\\
\end{array}\]
Using (\ref{Sab}) we see that the two columns are the same when $\x_a=\x_b$.
\sp

\noi{\bf The case \boldmath$(-,\,+)$}: The only $S$-factor involving only $\pm a$ and $\pm b$ is $\SX{b}{a}$ for $\s'$, which equals $-1$ when $\x_a=\x_b$. For $c\ne\pm a,\,\pm b$ we have the table
\[\begin{array}{lllll}
c<-b && \SX{a}{c}\,\SX{-a}{c}\,\SX{-c}{b} && \SX{b}{c}\,\SX{-b}{c}\,\SX{-c}{a} \\
-b<c<-a && \SX{-a}{c}\,\SX{a}{c}&& \SX{-c}{a}\,\SX{b}{c} \\
-a<c<a && \SX{a}{c}&& \SX{b}{c}\\
a<c<b && 1 && \SX{b}{c}\,\SX{c}{a} \\
c>b  &&  \SX{c}{b}&& \SX{c}{a}\\
\end{array}\]
Using (\ref{Sab}) and (\ref{Sprod}) we see that the two columns are the same when $\x_a=\x_b$.

To recapitulate, each nonzero $I(\s)$ with $\s\in\BSN$ is a sum of lower-order integrals in each of which some $\x_a=\x_b$. Denote by $I_{(a,\,b)}(\s)$ the sum of the lower-order integrals for $\s$, if any,  in which $\x_a=\x_b$. It follows from what we have shown that if we pair $\s$ and $\s'$ if they agree except for the positions of $\pm a$ and $\pm b$, and $a$ and $b$ are interchanged, then $I_{(a,\,b)}(\s)+I_{(a,\,b)}(\s')=0$. Therefore 
\[\sum_{\s}I_{(a,\,b)}(\s)=0\]
for each $(a,\,b)$. Summing over all pairs $(a,\,b)$ gives (\ref{sum}). This completes the proof of the theorem.
\sp 

We now expand on the discussion following (\ref{example}), and show that after the variable change the domains of integration for $\s$ and $\s'$ become the same.
The initial domain of integration (\ref{domain}), after the substitutions $\x_a\to\x_a+1/2q$, we write as
\[(\x_a)_{a\le N}\in\bigcup_{\m\in\S_N}\ \prod_{a\le N}\C_{\m(a)}.\]
Now the $\C_a$ are circles with center zero. Because of the way the $R_a$ were ordered, when we integrate with respect to $\abxs{j_0}$ by expanding its contours we pass a pole at 
$-\abxs{k_1}$ (where $k_1$ equals an $i$ or $j$) only for those contours for which 
$\m(|\s|(k_1))>\m(|\s|(j_0))$. This is a condition on $\m$, and our new domain of integration is a union of contours over only those $\m$ satisfying this condition:
\[(\x_{|\s|(\l)})_{\l\ne j_0}\in \bigcup_{\m(|\s|(k_1))>\m(|\s|(j_0))}\ \prod_{\l\ne j_0}\C_{\m(|\s|(\l))}.\]
Then we integrate with respect to $\abxs{k_1}$, and pass a pole at $\pm\abxs{k_2}$ only for those contours for which $\m(|\s|(k_2))>\m(|\s|(k_1))$. The new domain of integration is a union over fewer $\m$: 
\be(\x_{|\s|(\l)})_{\l\ne j_0,\,k_1}\in \bigcup_{{\m(|\s|(k_1))>\m(|\s|(j_0))\atop \m(|\s|(k_2))>\m(|\s|(k_1))}}\ \prod_{\l\ne j_0,\,k_1}\C_{\m(|\s|(\l))}.\label{domain1}\ee

If we do the same for $\s'$ and take the the same $k_1$ and $k_2$, the corresponding domain of integration would be a union over different $\m$ of different contours:
\[(\x_{|\s'|(\l)})_{\l\ne j_0,\,k_1}\in \bigcup_{{\m(|\s'|(k_1))>\m(|\s'|(j_0))\atop \m(|\s'|(k_2))>\m(|\s'|(k_1))}}\ \prod_{\l\ne j_0,\,k_1}\C_{\m(|\s'|(\l))}.\]

Suppose $\s$ and $\s'$ are $(a,\,b)$-paired. Switching the variables $\x_a$ and $\x_b$ in $\s'$ has the effect of replacing this by
\[(\x_{|\s|(\l)})_{\l\ne j_0,\,k_1}\in \bigcup_{{\m(|\s'|(k_1))>\m(|\s'|(j_0))\atop \m(|\s'|(k_2))>\m(|\s'|(k_1))}}\ \prod_{\l\ne j_0,\,k_1}\C_{\m(|\s'|(\l))}.\]
Let $\n$ be the permutation in $\S_N$ that interchanges $a$ and $b$ and leaves the rest of $[1,\,N]$ fixed. If we replace $\m$ by $\m\,\n$ on the right side (which we may do since $\m$ denoted a generic permutation) then we obtain precisely (\ref{domain1}). This is what we meant by the domains of integration for $\s$ and $\s'$ becoming the same after the variable switch.

\begin{center}{\bf Acknowledgment}\end{center}

The authors thank Neil O'Connell for valuable help in formulating the ASEP on the half-line.

This work was supported by the National Science Foundation through grants DMS-0906387 (first author) and DMS-0854934 (second author).

\end{document}